\baselineskip=14pt
\font\Bf=cmbx12

\def\X{X^*}
\def\tto{\,\scriptstyle$\lower1pt\vbox{\hbox{$\to$}\kern-10pt\hbox{$\to$}}$\,}
\def\co{{\rm co}}

\def\int{{\rm int }}
\def\co{{\rm co}}

\centerline {\Bf {Rockafellar's Sum Theorem}}

\bigskip

 Let $X$ be a Banach space,  $T, S:X\tto \X $ be maximal monotone operators. In [R] Rockafellar conjectured that $T+S$ is also maximal  provided that $D_T\cap\int(D_S) \ne \emptyset$ (here $D_T$ stands for the domain of $T$).  Theorem 3 in [EW] states that Rockafellar's conjecture is true provided that $D_T$ or $D_S$ is bounded.  In this note we show how to remove this restriction. 

\medskip
\noindent{\bf Theorem.}  Let $X$ be a Banach space,  $T, S:X\tto \X $ be maximal monotone operators. Assume that

\smallskip
(a) $D_S$ is convex and $\bigcup_{\lambda>0}\lambda(\co D_T-\co D_S)=X $ or (b) $D_T\cap\int(D_S) \ne \emptyset$. 

\smallskip\noindent
Then $ T + S$ is a maximal monotone operator.

\bigskip
\noindent {\it Proof.}  WLG we can assume that $0\in D_T \cap D_S$ if (a) is true or that  $0\in D_T \cap \int(D_S)$ if we (b) is true. Let $(x, x^*)$ be monotonically related to the graph of $T + S$. Choose a ball $B$ in $X$ (centered at $0$)  that contains $x$ and such that $ D_T \cap  D_S \cap B \ne \emptyset$. It is easily seen that 

\medskip
\noindent
(i) \qquad$(x, x^*)$ is related to $T + S + \partial I_B$. 

\medskip
\noindent
We shall now show that 

\smallskip\noindent
(ii) \qquad $\bigcup_{\lambda>0}\lambda(\co D_T-\co D_{S+\partial I_B})=X $.

\smallskip\noindent
This is obvious if (b) is true. So assume that (a) is true and let $z\in X$. Then there exist $\lambda>0$, $u\in \co D_T$ and $v\in \co D_S=D_S$ such that $z=\lambda(u-v)$.
Since $D_S$ is convex and $0\in D_T\cap D_S$, there exists $\mu$, $0<\mu<1$ such that $\mu v\in D_S\cap B$ and $\mu u\in \co D_T$. Then $z=(\lambda/\mu)(\mu u-\mu v)$ and thus (ii) is proved.

\medskip 
According to Theorem 3 mentioned above,  $ S + \partial I_B$ is maximal monotone. Because of (ii) the same theorem implies that $T + S + \partial I_B$ is maximal monotone, hence (because of (i)) $x \in  D_T \cap   D_S$. It is well known that this implies that  $T + S$ is a maximal monotone operator (see [V, Theorem 3.4 and Corollary 5.6] or [S, Theorem 24.1]).

\bigskip\noindent
{\bf References}

\medskip
\parindent=30pt

\item{[EW]} A. Eberhard \& R. Wenczel: All maximal monotone operators in a Banach space are of type FPV, Set Valued Var. Anal 22 (2014) 597-615. 

\item{[R]\ }  On the maximality of sums of nonlinear monotone operators, Trans. Amer. Math. Soc. 159 (1970), 81-99.

\item{[S]\ }  S. Simons: From Hahn-Banach to Monotonicity, Second Edition, Lecture Notes in Mathematics 1693 (2008), Springer-Verlag.

\item{[V]\ }  M. D. Voisei: The sum and chain rules for maximal monotone operators, Set-Valued Analysis, 16 (2008), 461--476.

\bigskip
\noindent {\bf Andrei $\&$ Maria Elena Verona}

\noindent 
verona@usc.edu

\bye